\documentclass[12pt,leqno]{amsart}
\usepackage{amssymb,verbatim,enumerate,ifthen,cite}
\usepackage[mathscr]{eucal}
\usepackage{cite}
\usepackage[utf8]{inputenc}
\usepackage[T1]{fontenc}
\def\N{\mathbb{N}}
\def\R{\mathbb{R}}

\def\M{\mathscr{M}}
\def\A{\mathscr{A}}

\newtheorem{theorem}{Theorem}[section]
\newtheorem*{theorem*}{Theorem}
\def\Thm#1#2{\ifthenelse{\equal{#1}{*}}{\begin{theorem*}#2\end{theorem*}}
	{\begin{theorem}\label{T#1}#2\end{theorem}}}

\newtheorem{proposition}[theorem]{Proposition}
\newtheorem*{proposition*}{Proposition}
\def\Prp#1#2{\ifthenelse{\equal{#1}{*}}{\begin{proposition*}#2\end{proposition*}}
	{\begin{proposition}\label{P#1}#2\end{proposition}}}
\def\prp#1{Proposition~\ref{P#1}}

\newtheorem{corollary}[theorem]{Corollary}
\newtheorem*{corollary*}{Corollary}
\def\Cor#1#2{\ifthenelse{\equal{#1}{*}}{\begin{corollary*}#2\end{corollary*}}
	{\begin{corollary}\label{C#1}#2\end{corollary}}}
\def\cor#1{Corollary~\ref{C#1}}

\newtheorem{lemma}[theorem]{Lemma}
\newtheorem*{lemma*}{Lemma}
\def\Lem#1#2{\ifthenelse{\equal{#1}{*}}{\begin{lemma*}#2\end{lemma*}}
	{\begin{lemma}\label{L#1}#2\end{lemma}}}

\theoremstyle{definition}
\newtheorem{remark}[theorem]{Remark}
\newtheorem*{remark*}{Remark}
\def\Rem#1#2{\ifthenelse{\equal{#1}{*}}{\begin{remark}\rm #2\end{remark}}
	{\begin{remark}\label{R#1}\rm #2\end{remark}}}

\newtheorem{example}[theorem]{Example}
\newtheorem*{example*}{Example}
\def\Exa#1#2{\ifthenelse{\equal{#1}{*}}{\begin{example*}\rm #2\end{example*}}
	{\begin{example}\label{Ex#1}\rm #2\end{example}}}

\def\eq#1{{\rm(\ref{E#1})}}
\def\Eq#1#2{\ifthenelse{\equal{#1}{*}}
	{\begin{equation*}\begin{aligned}[]#2\end{aligned}\end{equation*}}
	{\begin{equation}\begin{aligned}[]\label{E#1}#2\end{aligned}\end{equation}}}

\begin{document}
	\begin{flushright}
	\end{flushright}
	\vspace{5mm}
	
	\date{\today}
	
\title[Approximate monotonicity, subadditivity, and convexity in weighted topologies]
{Approximate monotonicity, subadditivity, and convexity in weighted topologies}
\author[A. R. Goswami]{Angshuman R. Goswami}
\address[A. R. Goswami]{Department of Mathematics, Faculty of Technical Informatics, University of Pannonia,
	H-8200 Veszprem, Hungary}
\email{goswami.angshuman.robin@mik.uni-pannon.hu}
\begin{abstract}
The central question of this paper is the following:\\
\newline
“If a topological space equipped with a weight function satisfies a certain property only approximately, can we recover the original structure without significantly altering the weights?”\\
\newline
To each open set of a topological space, we assign a non-negative weight and study approximate versions of three natural properties: monotonicity, subadditivity, and convexity. Through this study, we develop Hyers–Ulam-type stability results in general topology.  In addition, we investigate several topology-specific cases, as well as present minorant and sandwich-type results.

\end{abstract}

\maketitle

\section{Introduction}

In a broad sense, a Hyers-Ulam-type stability theorem can be expressed as follows:
\\

Let $\varepsilon$ denote the error with which a function $f$ with a certain property (say $P$) satisfies a given equation or inequality. If $\varepsilon$ is sufficiently small, then there exists an exact solution $g$ of the equation(inequality) such that the distance between the functions $f$ and $g$ is controlled by the size of the error $\varepsilon$. Symbolically, this can be represented as 
\Eq{*}{
\text{ equation/inequality with small error}
\quad\Longrightarrow\quad
\text{small distance to an exact solution}.
}
The origin of Hyers--Ulam-type stability goes back to a problem posed by Stanisław Ulam in 1940 concerning the stability of group homomorphisms. This question was answered by D. H. Hyers in 1941. Since then, researchers have extended this idea to various classes of functions satisfying properties such as monotonicity, convexity, subadditivity, periodicity etc., which are often expressed through some equations or inequalities. Some of the findings that inspire and impact our investigation are listed below in their simplified formulations.

\Thm{*}{\normalfont{\textbf{[Hyers (see \cite{Ulam})]}}
Let $\varepsilon\geq 0$ be fixed. Suppose $f:\R\to\mathbb{R}$ satisfies the following functional inequality
\Eq{*}{
|f(x+y)-f(x)-f(y)|\leq \varepsilon.
}
Then there exists a unique additive function $\A:I\to\mathbb{R}$ such that
$\|f-\A\|\leq \varepsilon/2$ holds.
}

\Thm{*}{\normalfont{\textbf{[Hyers, Ulam (see \cite{Hyers})]}} We assume $\varepsilon>0$ be fixed and $f:\R\to\R$ such that for every $x,y\in \R$ and for all $[0,1]$, the following functional inequality is satisfied
\Eq{*}{
f(tx+(1-t)y)\leq tf(x)+(1-t)f(y)+\varepsilon.
}
Then there exists a convex function $C:\R\to\R$ such that 
$\|f-C\|\leq \varepsilon/2$ holds.
}

\Thm{*}{\normalfont{\textbf{[P\'ales(see \cite{Pales})]}}
Let $\varepsilon>0$ be fixed. If the function $f:\R\to\R$ satisfies the following functional inequality
\Eq{*}{
f(x)\leq f(y)+\varepsilon \qquad \mbox{for all}\quad x,y\in \R\,,
}
then there exists a monotone(increasing) function $\M:\R\to\R$ such that 
$\|f-\M\|\leq \varepsilon/2$ holds.
}
\Thm{*}{\normalfont{\textbf{(see \cite{Goswami})}}
Let $\varepsilon>0$ be fixed and 
$f:\R_+\to\R_+$ such that for every $x\in \R_+$, the following inequality holds
\Eq{*}{
f(x)\leq f(x_1)+\cdots+f(x_n)+\varepsilon,\qquad \mbox{ where}\qquad \sum_{i=1}^nx_i=x.
}
Then, corresponding to $f$, there exists a subadditive minorant $g:\R_+\to\R_+$  that satisfies the inequality $0\leq f-g\leq \varepsilon $.
}

In addition, converse statements were established for most of the results mentioned above. Further background on these function classes and related topics can be found in the classical monographs \cite{Hardy,Kuczma}.\\

Over the last few decades, Hyers--Ulam-type stability has also been studied for ordinary and partial differential equations, difference equations, and complex functional inequalities. More recently, stability theorems have been developed in discrete settings, particularly for sequences.
The results obtained in our recent studies of Hyers--Ulam-type stability on graphs (see \cite{Ang}) naturally led us to consider a broader setting. The aim of this research note is to  formulate them in a more general and abstract framework.\\

Throughout this paper, the symbols $\N$ and $\R$ denote the sets of natural and non-negative real numbers, respectively. The notation $(X,\mathcal{T})$ denotes a finite topological space. The function $w:\mathcal{T}\to\R_+$ with $w(\varnothing)=0$ will be called weight function. The topological space $(X,\mathcal{T})$ together with the weight function $w$ is said to be a weighted topology and will be represented  by $(X,\mathcal{T},w)$. It is important to note that despite having the same topological structure, the two weighted topologies $(X,\mathcal{T},w_1)$ and $(X,\mathcal{T},w_2)$ are distinct in the sense that the assigned weight functions differ. For $\mathcal{T}_0\subset \mathcal{T}$, the notation $(X,\mathcal{T}_0,w_0)$ represents the weighted topology, where 
$w_0:\mathcal{T}_0\to\R$ is nothing but $w$ restricted to $\mathcal{T}_0$. Similarly, for any proper open set $Y\in\mathcal{T}$, the weight function $w_{_{Y}}:\mathcal{T}_{_{Y}}\to\R$ in the resultant relative topology $(Y,\mathcal{T}_{_{Y}},w_{_{Y}} )$ holds the same meaning of restrictive weights. Clearly, specifically allocated weight distribution enables the investigation of several fundamental structures within the topology
$\mathcal{T}$.  Prominent among these are monotonicity, subadditivity, and convexity. \\

For a fixed $\varepsilon>0$, a weighted topology $(X,\mathcal{T},w)$ is called \textit{approximately monotone}, if all
$G,G'\in\mathcal{T}$ satisfy the following discrete functional inequality

\Eq{0}{
	w(G)\leq w(G')+\varepsilon \qquad\mbox{provided}\qquad G\subset G'\,.
}
If $\varepsilon=0$, we call it a \textit{monotone (increasing) topology}. By definition, in an increasing topology for any $G\in\mathcal{T}$, the inequality $w(\varnothing)\leq w(G)\leq w(X)$ is evident. Besides, every discrete topology $(X,\mathcal{T},w)$ that satisfies the condition $w(G)+w(G^c)=1$ for all $G\in\mathcal{T}$, naturally induces a probability space. To illustrate the notion of approximately monotone topology, consider the topology $\mathcal{T}=\{\varnothing,
\{a\},\{b\},\{a,b\},X\}$ constructed on the set $X=\{a,b,c\}$. Under the associated weights $w(\varnothing)=0,
 w(\{a\})=0.5,$ $w(\{b\})=0.75, w(\{a,b\})=1.01$ and $w(X)=1$ it exhibits approximate monotonicity with 
$\varepsilon=0.01$. We demonstrate that every approximately monotone topology admits an ordinary monotone topology on the same underlying structure, with the associated weight functions differing by a minimal amount. A conversion assertion for this statement is also presented. Moreover, some results concerning finer, indiscrete, relative, partitioning topologies are included. \\

Let $\varepsilon>0$ be fixed. A weighted topology
$(X,\mathcal{T},w)$ is said to be \textit{approximately subadditive}, if for any  $G_1,\cdots ,G_n\in\mathcal{T}$, the weight function $w$ satisfies the following discrete inequality

\Eq{210}{
w(G)\leq w(G_1)+\cdots+w(G_n)+\varepsilon \qquad\mbox{where}\qquad G=\overset{n}{\underset{i=1}\cup} G_i\,.
}

Clearly, if $\varepsilon=0$, the weighted topology  $(X,{\mathcal{T}},w)$ is nothing but a \textit{subadditive} topology. Our result shows that any approximately subadditive topology can be replaced by a proper subadditive topology on the same space, having the same open sets, while keeping the distance between them bounded above by ${\varepsilon}/{2}$. Moreover, we presented several findings showing relationships between monotone and subadditive topologies.\\

On the other hand, for a fixed $\varepsilon>0$, a weighted topology $(X,\mathcal{T},w)$ is termed \textit{approximately convex}, if for any triplet of non-empty open sets ${G_*}, G, {G^*}\in\mathcal{T}$ with 
${G_*}\subset G \subset {G^*}$, the following discrete functional inequality holds
\Eq{99}{
2w(G)\leq w({G_*})+w({G^*})+\varepsilon\,\,.
}
In case $\varepsilon=0$, we say it is a \textit{convex topology}. This definition of convexity is inspired by the notion of convex sequences or discrete convex functions, which was introduced by Mitrinovi'c in his book \cite{Mitrinovic}. A salient feature of this topological convexity is that, for every chain
$
\mathcal{C}:=\{G_i,\subseteq\}
$
in the topology $(X,\mathcal{T})$, either the sequence \(\bigl(w(G_i)\bigr)\) is monotone, or there exists a proper open set \(G\in\mathcal{C}\) such that the corresponding weighted sequence from \(\varnothing\) to \(G\) 
(excluding 
$\varnothing$) is decreasing, while from \(G\) to \(X\) is increasing. We presented a result showing that if a weighted topology  $(X,\mathcal{T},w)$ exhibits approximate convexity, then there exists a non-trivial (non-constant)convex topology 
$(X,\mathcal{T},\overline{w})$ such that for all $G\in\mathcal{T}\setminus\{\varnothing\}$, the inequality below satisfies
\Eq{*}{
\min_{G\in\mathcal{T}\setminus\{\varnothing\}}\Bigg\{\max\bigg\{ w(G)-\dfrac{\varepsilon}{2}\,,0\,\bigg\}\Bigg\}\leq \overline{w}(G)\leq w(G)+\dfrac{\varepsilon}{2}.
}
Although, the above result not fall under the category of Hyers-Ulam-type stability; yet it indicates the possibility of local stability cantering around the open set that carries the least weight.
\\

To keep the presentation simple, this research note focuses exclusively on finite topologies. We expect that several of the results can be extended to suitable infinite topologies after appropriate technical modifications, such as replacing minima by infima and maxima by suprema in the relevant constructions. A complete treatment of this infinite-dimensional setting is left for future work. We start our investigation with an approximately monotone topology.
\section{Main Results}
First, we go through some basic structural results of monotone topologies. The proof of the assertions are straightforward, and hence we propose only the statements.
\Prp{010}{ Let $(X,\mathcal{T},w)$ be a  monotone topology. Then the following assertions hold:
\begin{enumerate}
\item If $(X,\mathcal{T})$ is a indiscrete topology, then $(X,\mathcal{T},w)$ is a monotone topology for all possible weight functions $w$.
\item If ${\mathcal{T}}_0\subset\mathcal{T}$, then the weighted topology $(X,{\mathcal{T}}_0,w_0)$ possesses the same monotonicity.
\item If $Y\in \mathcal{T}$, then the weighted relative topology $(Y,\mathcal{T}_{_{Y}},w_{_{Y}})$ possesses the same monotonicity. 
\item If $(X,\mathcal{T})$ is a partition topology, then $(X,\mathcal{T}_{_{P}},w')$ possesses the same monotonicity;
\newline
where $\mathcal{T}_{_{P}}:=\{G:X\setminus G\in\mathcal{T}\}$ and $w':\mathcal{T}_{_{P}}\to\R_+$ is $w'(X-G)=w(X)-w(G).$
\end{enumerate}
}
\begin{proof}
The proofs of these assertions are left to the reader.
\end{proof}
In the next proposition, we demonstrate how to formulate the largest monotone topological minorant for any weighted topological space.
\Prp{1} {If $(X,\mathcal{T},w)$ be a weighted topology, then
$(X,\mathcal{T},\widetilde{w})$ is the largest monotonically increasing minorant of $(X,\mathcal{T},w)$, where function $\widetilde{w}:\mathcal{T}\to\R_+$ is defined as follows:
\Eq{1}{
\widetilde{w}\big(G\big):=\min_{G'\in\mathcal{T}}\Big\{w\big(G'\big)\,\,\Big|\,\,G\subseteq G'\Big\}.
}
}
\begin{proof}
By definition, the following inequality is evident. 
\Eq{*}{
\widetilde{w}\big(G\big)=\min_{G'\in\mathcal{T}}\Big\{w\big(G'\big)\,\,\Big|\,\,G\subseteq G'\Big\}\leq w(G);
}
which shows that $(X,\mathcal{T},\widetilde{w})$ is a minorant of $(X,\mathcal{T},{w})$.\\

To establish the increasingness, we assume any two arbitrary open sets $G_1,G_2\in\mathcal{T}$ such that $G_1\subset G_2$ holds. Then using \eq{1}, we can conclude the following inequality
\Eq{*}{
\widetilde{w}\big(G_1\big)= \min_{G'\in\mathcal{T}}\Big\{w\big(G'\big)\,\,\Big|\,\,G_1\subseteq G'\Big\} \leq \min_{G'\in\mathcal{T}}\Big\{w\big(G'\big)\,\,\Big|\,\,G_2\subseteq G'\Big\}=\widetilde{w}\big(G_2\big).
}
This shows that $(X,\mathcal{T},\widetilde{w})$ is an increasing topology. Thus we are only required to show that no other weighted monotone(increasing) topology can lie in between  $(X,\mathcal{T},w)$ and $(X,\mathcal{T},\widetilde{w})$.\\

If possible, we assume that there exists another mapping $w_0:\mathcal{T}\to\R_+$ such that $(X,\mathcal{T},w_0)$ is the largest increasing minorant of $(X,\mathcal{T},w)$. This together with the monotonicity of $(X,\mathcal{T},\widetilde{w})$ yields that for any open set $G\in \mathcal{T}$, all the existing open sets $G'\in\mathcal{T}$ with $G\subseteq G'$, satisfy the following functional inequality
\Eq{*}{
\widetilde{w}\big(G\big)\leq{w_0}\big(G\big)\leq w_0\big(G'\big)\leq w\big(G'\big).
}
By taking the minimum to the right-most side of the above inequality, we have 
\Eq{*}{
{w_0}\big(G\big)\leq \min_{G'\in\mathcal{T}}\Big\{w\big(G'\big)\,\,\Big|\,\,\mbox{for all}\,\,G\subseteq G'\Big\}=\widetilde{w}\big(G\big).
}
This contradicts our assumption.
Therefore, $(X,\mathcal{T},\widetilde{w})$ is the largest monotone(increasing) minorant of $(X,\mathcal{T},w).$ 
\end{proof}
In the following proposition, we show that if two weighted topologies exhibit some specific inequality conditions, then it is feasible to sandwich a monotone topology in between them.
\Cor{1}{Let $(X,\mathcal{T},{w_1})$ and $(X,\mathcal{T},{w_2})$ are two weighted topologies such that for every 
$G\in \mathcal{T}$, the inequality $w_1\big(G\big)\leq w_2\big(G'\big)$ holds for all $G'\in \mathcal{T}$ with $G\subseteq G'$. Then there exists a monotone (increasing) topology $(X,\mathcal{T},\widetilde{w})$ sandwiched in between  $(X,\mathcal{T},{w_1})$ and $(X,\mathcal{T},{w_2})$.
}
\begin{proof}
We can construct a monotone weight function $\widetilde{w}:\mathcal{T}\to\R_+$ as in \eq{1}, which ensures the following inequality
\Eq{20}{
w_1\big(G\big)\leq \widetilde{w}\big(G\big)\leq w_2\big(G\big)\qquad\mbox{for all}\qquad G\in \mathcal{T}.
} 
This validates the sandwich property and completes the proof.
\end{proof}
The following result can be interpreted as a Hyers-Ulam stability theorem for monotone topological spaces.
\Prp{2.2}{
Let $(X,\mathcal{T},w)$ be an approximately monotone topology, then there exists a monotone(increasing) topology $(X,{\mathcal{T}},\widetilde{w})$ such that $\|w-\widetilde{w}\|\leq \varepsilon/2$ holds. Conversely if $(X,{\mathcal{T}},\widetilde{w})$ is a monotone topology satisfying the inequality $\|w-\widetilde{w}\|\leq \varepsilon/2$, then $(X,\mathcal{T},w)$ is approximately monotone.
}
\begin{proof}
The inequality \eq{0} can also be re-phrased as follows
\Eq{10}{
	w(G)-\dfrac{\varepsilon}{2}\leq w(G')+\dfrac{\varepsilon}{2} \qquad\mbox{provided}\qquad G\subseteq G'\,.
}
To maintain simplicity, we introduce the following weight functions on the topology $\mathcal{T}$.
\Eq{*}{
w_1(G):=\max\bigg\{0,\,w(G)-\dfrac{\varepsilon}{2}\bigg\}\qquad\mbox{and}\qquad w_2(G):=w(G)+\dfrac{\varepsilon}{2}.
}
Thus the inequality \eq{10} can be restructured as follows
\Eq{*}{
w_1(G)\leq w_2(G') \qquad\mbox{provided}\qquad G\subseteq G'\,.
}
As a consequence of established inequality \eq{20} of \cor{1}, we obtain the following inequality 
\Eq{90}{
w(G)-\dfrac{\varepsilon}{2}\leq \widetilde{w}(G) \leq w(G)+\dfrac{\varepsilon}{2} \qquad \mbox{for all}\quad G\in \mathcal{T}.
}
This yields the norm inequality $\|w-\widetilde{w}\|\leq \varepsilon/2$ and completes the proof of first assertion.\\

To show the converse part, we assume that \eq{90} holds where 
$\widetilde{w}:\mathcal{T}\to\R_+$ is monotone(increasing). Now, we consider $G,G'\in \mathcal{T} $ with $G\subset G'.$ Then using \eq{90}, we can compute the following inequality
\Eq{*}{
w(G)\leq \widetilde{w}(G)+\dfrac{\varepsilon}{2}
\leq \widetilde{w}(G')+\dfrac{\varepsilon}{2}
\leq \bigg({w}(G')+\dfrac{\varepsilon}{2}\bigg)+\dfrac{\varepsilon}{2}={w}(G')+\varepsilon.
}
This shows the approximate monotonicity of $(X,\mathcal{T},w)$ and completes the proof.
\end{proof}
Next, we propose a stability result for subadditive topologies. For clarity, rather than presenting a lengthy proof, we decompose our findings into a sequence of smaller, intermediate results following a traditional framework. First, we construct the largest subadditive minorant for any weighted topological space. Second, we establish the underlying conditions under which a proper subadditive topology can act as a separator for two weighted topological spaces. Finally, utilizing these components, we establish our Hyers-Ulam stability result for subadditive topologies.

\Prp{2}{ If $(X,\mathcal{T},w)$ is a weighted topology, then  $(X,\mathcal{T},\widehat{w})$ is the largest subadditive minorant of $(X,\mathcal{T},w)$, where the function $\widehat{w}:\mathcal{T}\to\R_+$ is defined as follows:
\Eq{3}{
\widehat{w}\big(G\big):=\min\Big\{w\big(G_1\big)+\cdots +w\big(G_n\big)\,\,\Big|\,\,G_1,\cdots,G_n\in \mathcal{T}\,\,\mbox{such that}\,\,\overset{n}{\underset{i=1}{\cup}}G_i=G \}.
}
Additionally, if $(X,\mathcal{T},w)$ is monotone then $(X,\mathcal{T},\widehat{w})$ also possesses same monotonicity.
}
\begin{proof}
The minorant part is clearly evident from the following inequality
\Eq{*}{
\widehat{w}\big(G\big)=\min\Big\{w\big(G_1\big)+\cdots +w\big(G_n\big)\,\,\Big|\,\,G_1,\cdots,G_n\in\mathcal{T}\,\,\mbox{such that}\,\,\overset{n}{\underset{i=1}{\cup}}G_i=G \}\leq w\big(G\big).}

Now to establish the subadditivity of 
$(X,\mathcal{T},\widehat{w})$, we start by assuming two distinct arbitrary open sets $G_1$ and $G_2$ of 
$\mathcal{T}$, which definitely implies $G_1\cup G_2\in \mathcal{T}$. Then from \eq{3}, we can obtain the following two 
\Eq{4}{
\widehat{w}(G_1)=w\big(G^{^{1}}_{_{1}}\big)+\cdots +w\big(G^{^{1}}_{_{k}}\big)\,\,&\mbox{such that}\,\,G^{^{1}}_{_{1}},\cdots,G^{^{1}}_{_{k}}\in \mathcal{T}\,\,\mbox{satisfying}\,\,\overset{k}{\underset{i=1}{\cup}}G^{^{1}}_{_{i}}=G_1\\
&\quad\,\,\mbox{and}\\
\widehat{w}(G_2)=w\big(G^{^{2}}_{_{1}}\big)+\cdots +w\big(G^{^{2}}_{_{\ell}}\big)\,\,&\mbox{such that}\,\,G^{^{2}}_{_{1}},\cdots,G^{^{2}}_{_{\ell}}\in \mathcal{T}\,\,\mbox{satisfying}\,\,\overset{\ell}{\underset{i=1}{\cup}}G^{^{2}}_{_{i}}=G_2.
}
Also, one can easily observe the following equality
\Eq{*}{G^{^{1}}_{_{1}}\cup\cdots\cup G^{^{1}}_{_{k}}\cup G^{^{2}}_{_{1}}\cup\cdots\cup G^{^{2}}_{_{\ell}}=G_1\cup G_2\in \mathcal{T}.
}
In the framework of \eq{3}, we can conclude the following
\Eq{*}{
\widehat{w}\big(G_1\cup G_2\big)\leq w\big(G^{^{1}}_{_{1}}\big)+\cdots +w\big(G^{^{1}}_{_{k}}\big)+\cdots+w\big(G^{^{2}}_{_{1}}\big)+\cdots +w\big(G^{^{2}}_{_{\ell}}\big).
}
Using \eq{4} in the above inequality, we get 
\Eq{*}{\widehat{w}\big(G_1\cup G_2\big)\leq w\big(G_1\big)+w\big(G_2\big).}
This shows $(X,\mathcal{T},\widehat{w})$ possesses subadditivity.\\

Now, we assume that there exist another mapping $w_0:\mathcal{T}\to\R_+$ such that $(X,\mathcal{T},{w_0})$ is the largest subadditive minorant of $(X,\mathcal{T},{w})$. Hence, for arbitrarily chosen $G\in \mathcal{T}$ together with $G_1,\cdots,G_n\in\mathcal{T}$ satisfying $\overset{n}{\underset{i=1}{\cup}}G_i=G$, we can obtain the following discrete functional inequalities
\Eq{*}{
w_0\big(G\big)\leq w_0\big(G_1\big)+\cdots +w_0\big(G_n\big)\leq w\big(G_1\big)+\cdots +w\big(G_n\big).
}
Taking the minimum to the right-most part of the above inequality, we have
\Eq{*}{
w_0\big(G\big)
\leq\min\Big\{ w\big(G_1\big)+\cdots +w\big(G_n\big)\,\,\Big|\,\,G_1,\cdots,G_n \in \mathcal{T}\,\,\mbox{satisfying}\,\,\overset{n}{\underset{i=1}{\cup}}G_i=G\Big\}=\widehat{w}\big(G\big).
}
This contradicts our assumption on $(X,\mathcal{T},{w_0})$. Hence, $(X,\mathcal{T},\widehat{w})$ is the largest subadditive minorant of $(X,\mathcal{T},w)$. This establishes the first part of the proposition.\\

To validate the second part of the assertion, without loss of generality, we consider that $(X,\mathcal{T},{w})$ is monotonically increasing. Let $G_1,G_2\in \mathcal{T}$ such that $G_1\subset G_2$. Then from \eq{3} of \prp{2}, there must exists open sets $G_1^{1},\cdots,G_k^{1}$ $\subseteq G_1$ and $G_1^{2},\cdots,G_\ell^{2}$ $\subseteq G_2$ of $(X,\mathcal{T})$ such that both equalities of \eq{4} holds along with the mentioned conditions. Also, one can easily observe the following inclusions
\Eq{*}{
G_1^{2}\cap G_1\subset G_1,\cdots, G_\ell^{2}\cap G_1\subset G_1 \quad\mbox{which also implies}\quad \overset{\ell}{\underset{i=1}{\cup}}\big(G_i^{2}\cap G_1\big)=G_1.
}
Using these and then utilizing the monotonicity of $(X,\mathcal{T},{w})$, the first equality of \eq{4} can be extended as follows
\begin{small}
\Eq{*}{
\widehat{w}(G_1)=w\big(G_1^{1}\big)+\cdots +w\big(G_k^{1}\big)
\leq w\big(G_1^{2}\cap G_1\big)+\cdots +w\big(G_\ell^{2}\cap G_1\big)
\leq  w\big(G_1^{2}\big)+\cdots +w\big(G_\ell^{2}\big)
=\widehat{w}(G_2).
}
\end{small}
Since $G_1,G_2\in \mathcal{T}$ are arbitrary, the above inequality establishes increasingness of $(X,\mathcal{T},\widehat{w})$.\\

Due to our assumption on $w$, the weighted topology $(X,\mathcal{T},w)$ can not be monotonically decreasing. However, we can consider $\mathcal{T}\setminus\{\varnothing\}$ monotonically decreasing and proceed by considering two $G_1,G_2\in\mathcal{T}\setminus\{\varnothing\}$ such that $G_1\subset G_2$. This implies the following two inequalities
\Eq{*}{
w\big(G_1\cup G_2\big)\leq w\big(G_1\big)\qquad \mbox{and}\qquad w\big(G_1\cup G_2\big)\leq w\big(G_2\big).
}
Thus $w\big(G_1\cup G_2\big)\leq w\big(G_1\big)+w\big(G_2\big)$ is obvious. That is, $(X,\mathcal{T},{w})$ possesses subadditivity. In other words,
if $(X,\mathcal{T},w)$ is decreasing in $\mathcal{T}\setminus\{\varnothing\}$, then $(X,\mathcal{T},{w})=(X,\mathcal{T},\widehat{w})$.\\

This validates the second assertion and completes the proof.
\end{proof}
The following corollary is a direct consequence of \prp{2}. Therefore, we provide only the draft of the proof.
\Cor{2}{If $(X,\mathcal{T},w_1)$ and $(X,\mathcal{T},w_2)$ are two weighted topologies such that 
\Eq{*}{w_1\big(G\big)\leq w_2\big(G_1\big)+\cdots+w_2\big(G_n\big)
}
holds for any $G_1,\cdots,G_n\in \mathcal{T}$ with $\overset{n}{\underset{i=1}{\cup}}G_i=G$,
then there exists a subadditive topology $(X,\mathcal{T},\widehat{w})$, that satisfies the following discrete functional inequality
\Eq{*}{
w_1\big(G\big)\leq \widehat{w}\big(G\big)\leq w_2\big(G\big)\qquad\mbox{for all}\qquad G\in \mathcal{T}.
}
}
\begin{proof}
We construct the weight function $\widehat{w}:\mathcal{T}\to\R_+$ as in \eq{3}. From there, the inequality $w_1\big(G\big)\leq \widehat{w}\big(G\big)\leq w_2\big(G\big)$, for all $G\in \mathcal{T}$ is obvious. The proof of the remaining assertions related to subadditivity are analogous to the first part of \prp{2}.
\end{proof}
The next result can be treated as the Hyers-Ulam-type stability for subadditive topology.
\Cor{107}{
Let $(X,\mathcal{T},w)$ be an approximately subadditive topology. Then there exists a subadditive topology 
$(X,\mathcal{T},\widehat{w})$ such that the norm inequality
$\|w-\widehat{w}\|<{\varepsilon}/2$ holds.}

\begin{proof}
Since, $(X,\mathcal{T},w)$ is an approximately subadditive topology, we first rewrite the inequality \eq{210} along with the mentioned conditions as follows 
\Eq{*}{
w(G)-\dfrac{\varepsilon}{2}\leq w(G_1)+\cdots+w(G_n)+\dfrac{\varepsilon}{2} \qquad\mbox{where}\qquad G=\overset{n}{\underset{i=1}\cup} G_i\,.
}
Now, in the above inequality taking the minimum with respect to all $G_{i}'s\in \mathcal{T}$ satisfying the condition $G=\overset{n}{\underset{i=1}\cup} G_i$ we arrive at the following inequality
\Eq{*}{
w(G)-\dfrac{\varepsilon}{2}\leq \min\Big\{ w(G_1)+\cdots+w(G_n)+\dfrac{\varepsilon}{2}\Big\}=\min\Big\{ w(G_1)+\cdots+w(G_n)\Big\}+\dfrac{\varepsilon}{2}\leq w(G)+\dfrac{\varepsilon}{2} \,.
}
We know that, algebraic summation of a non-negative subadditive function and a constant again results in a subadditive function. Thus in view of \cor{2}, we can have the subadditive weight 
$\widehat{w}:\mathcal{T}\to\R_+$, setting $\widehat{w}(\varnothing)=0$, that satisfies the following inequality
\Eq{*}{
w(G)-\dfrac{\varepsilon}{2}\leq \widehat{w}(G)\leq w(G)+\dfrac{\varepsilon}{2} \qquad \qquad({G\in \mathcal{T}}) \,.
}
This shows the validity of the assertion.
\end{proof}
Although the next result is not a Hyers--Ulam-type stability theorem in the classical sense, nonetheless it demonstrates that an approximately convex topology can still be associated with an exact convex topology satisfying suitable quantitative bounds. Most importantly, both the upper and lower bounds are governed by the nature of the approximately convex topology.

\Prp{4812}
{If the weighted topology $(X,\mathcal{T},w)$ is approximately convex, then there exists a non-trivial (non-constant) convex topology $(X,\mathcal{T},\overline{w})$ such that following inequality satisfies
\Eq{369}{
\min_{G\in\mathcal{T}\setminus\{\varnothing\}}\Bigg\{\max\bigg\{ w(G)-\dfrac{\varepsilon}{2}\,,0\,\bigg\}\Bigg\}\leq \overline{w}(G)\leq w(G)+\dfrac{\varepsilon}{2} \qquad \Big( G\in\mathcal{T}\setminus\{\varnothing\}\Big).
}
}
\begin{proof}
We will establish the proposition in several steps.
First, we prove the following statement:\\

"For any weighted topology $(X,\mathcal{T},w)$, there exists a non-trivial (non-constant) convex topology $(X,\mathcal{T},\overline{w})$ such that the inequality $ \overline{w}(G)\leq w(G)$ holds for all $G\in\mathcal{T}.$"\\

To validate it, let  $G\in\mathcal{T}$ be arbitrary. We construct a non-negative sequence $\Big(w_n(G)\Big)_{n=0}^{\infty}$ as follows
\Eq{119}{
w_0(G):=\min_{{G_*}\subseteq G\subseteq{G^*}} \dfrac{w({G_*})+w({G^*})}{2}\quad\mbox{and}\quad  w_n(G):=\min_{{G_*}\subseteq G\subseteq{G^*}} \dfrac{w_{n-1}({G_*})+w_{n-1}({G^*})}{2}\quad (n\in\N).
}
By definition, it can be easily observed that if $G=\varnothing$, then $\overline{w}(\varnothing)={w}(\varnothing)=0.$
One can also verify that the sequence $\Big(w_n(G)\Big)_{n=0}^{\infty}$ is monotone(decreasing) and bounded below by $0$. This ensures the convergence of the sequence. Now, we define the weight function $\overline{w}:\mathcal{T}\to\R_+$ as
\Eq{*}{
  \overline{w}(G)=\lim_{n\to\infty} w_n(G).
}
Hence $\overline{w}\leq w$ is evident. We need to show that the weight function $\overline{w}$ possesses convexity.\\

We assume $G\in\mathcal{T}\setminus\{\varnothing,X\}$ be fixed, and $G_1, G_2\in\mathcal{T}\setminus\{\varnothing\}$ are two arbitrarily selected open sets such that $G_1\subset G\subset G_2$ holds. Then by definition of $\overline{w}$, we can compute the following 

\Eq{*}{
 \overline{w}(G)
= \lim_{n\to\infty} w_n(G)
&= \lim_{n\to\infty}\Bigg(\min_{{G_*}\subseteq G\subseteq{G^*}} \dfrac{w_{n-1}({G_*})+w_{n-1}({G^*})}{2}\Bigg)\\
&\leq\lim_{n\to\infty}\Bigg(\frac{w_{n-1}(G_1)+w_{n-1}(G_2)}{2}\Bigg)=\dfrac{\overline{w}(G_1)+\overline{w}(G_2)}{2}.
}
This shows that $(X,\mathcal{T},\overline{w})$ is a convex minorant of the weighted topology $(X,\mathcal{T},{w})$ and establishes the statement.
\\

Now, we have the most required tool to establish the theorem. We start by re-structuring the inequality \eq{99} as below
\Eq{*}{
\max\bigg\{w(G)-\dfrac{\varepsilon}{2}\,,\,0\bigg\}
\leq \underset{G_*,G*\in\mathcal{T}\setminus\{\varnothing\}}{\min_{{G_*}\subseteq G\subseteq{G^*}}}\dfrac{\Big(w({G_*})+\varepsilon/2\Big)+\Big(w({G^*})+\varepsilon/2\Big)}{2}\,\qquad \Big({G\in\mathcal{T}\setminus\{\varnothing\}}\Big).
}
For simplicity, we use substitutions to represent the above inequality as follows
\Eq{130}{
w'(G)\leq w''(G)\qquad\mbox{for all}\qquad G\in\mathcal{T}\setminus\{\varnothing\}.
}\\
where the weight functions $w', w'':\mathcal{T}\setminus\{\varnothing\}\to\R_+$ are given by 
\begin{small}
\Eq{212}{
w'(G):=\max\bigg\{w(G)-\dfrac{\varepsilon}{2}\,,\,0\bigg\}\quad\mbox{and}\quad w''(G):=\underset{G_*,G*\in\mathcal{T}\setminus\{\varnothing\}}{\min_{{G_*}\subseteq G\subseteq{G^*}}}\dfrac{\Big(w({G_*})+\varepsilon/2\Big)+\Big(w({G^*})+\varepsilon/2\Big)}{2}\,.
}
\end{small}
Now, for each $G\in\mathcal{T}\setminus\{\varnothing\}$ we define a sequence $\Big(w''_n(G)\Big)_{n=0}^{\infty}$ (analogous to \eq{119}) as follows
\Eq{1119}{
w''_0(G):=w''(G)\qquad\mbox{and}\qquad  w''_n(G):=\underset{G_*,G*\in\mathcal{T}\setminus\{\varnothing\}}{\min_{{G_*}\subseteq G\subseteq{G^*}}} \dfrac{w''_{n-1}({G_*})+w''_{n-1}({G^*})}{2}\,.
}
Thus following similar mathematical arguments, we can show the existence of a non-trivial convex weight function $\overline{w}:\mathcal{T}\to\R_+$ such that for all open sets $G\in\mathcal{T}\setminus\{\varnothing\}$, the following functional inequality satisfied
\Eq{234}{
\overline{w}(G)\leq w''_{n}(G)\leq w''(G)\qquad (n\in\N).
}
Also from the inequality \eq{130} and \eq{234}, for all $G\in\mathcal{T}\setminus\{\varnothing\}$, we conclude the following inequality
\Eq{120}{
\min_{G\in\mathcal{T}\setminus\{\varnothing\}}w'(G)\leq \overline{w}(G).
}
If not, then there exists a $G\in\mathcal{T}\setminus\{\varnothing\}$ such that $\overline{w}(G)<\underset{G\in\mathcal{T}\setminus\{\varnothing\}}{\min}w'(G)$ holds. In other words, for all $G\in\mathcal{T}\setminus\{\varnothing\}$, this implies the validity of the following inequality
\Eq{*}{
\lim_{n\to\infty} w_n''(G)
= \lim_{n\to\infty}\Bigg(\,\,\underset{G_*,G*\in\mathcal{T}\setminus\{\varnothing\}}{\min_{{G_*}\subseteq G\subseteq{G^*}}} \dfrac{w_{n-1}''({G_*})+w_{n-1}''({G^*})}{2}\Bigg)<\underset{G\in\mathcal{T}\setminus\{\varnothing\}}{\min}w'(G).
}
In view of \eq{212} and \eq{1119} (besides $\mathcal{T}$ is a finite topology), for all $G\in\mathcal{T}\setminus\{\varnothing\}$ there exists an $n\in\N$ such that the following  inequality holds
\Eq{*}{
w''_n(G)=\underset{G_*,G*\in\mathcal{T}\setminus\{\varnothing\}}{\min_{{G_*}\subseteq G\subseteq{G^*}}} \dfrac{w''_{n-1}({G_*})+w''_{n-1}({G^*})}{2}\leq\min_{G\in\mathcal{T}\setminus\{\varnothing\}}w'(G).
}
This yields the existence of either an ${G_*}\subseteq G$ or an ${G^*}\supseteq G$ (the open sets $G_*,G*$ are non-empty) such that at least one of the following inequalities is satisfied
\Eq{*}{
w''_{n-1}\big({G_*}\big)\leq\min_{G\in\mathcal{T}\setminus\{\varnothing\}}w'(G)\qquad\mbox{or}\qquad
w''_{n-1}\big({G^*}\big)\leq\min_{G\in\mathcal{T}\setminus\{\varnothing\}}w'(G).
}
Without loss of generality, we assume that out of these two inequalities, the first inequality, that is $w''_{n-1}\big({G_*}\big)\leq\underset{G\in\mathcal{T}\setminus\{\varnothing\}}\min w'(G)$ holds along with the assumed conditions. Now using the same logic, we proceed, and eventually we get the existence of an $G_0\in\mathcal{T}\setminus\{\varnothing\}$ that satisfies the following inequality 
\Eq{*}{
w_0''(G_0)=w''(G_0)\leq\min_{G\in\mathcal{T}\setminus\{\varnothing\}}w'(G).
}
This results in a contradiction to the inequality \eq{130}. Hence, \eq{120} is valid.\\

By combining the two inequalities \eq{234} and \eq{120}, we arrive at 
\Eq{700}{
\min_{G\in\mathcal{T}\setminus\{\varnothing\}}w'(G)\leq \overline{w}(G)\leq w''(G).
}
Also from \eq{212}, we obtain the following two extensions
\Eq{*}{
&\qquad\quad\qquad\qquad
\min_{G\in\mathcal{T}\setminus\{\varnothing\}}w'(G)
=\min_{G\in\mathcal{T}\setminus\{\varnothing\}}\Bigg\{\max\bigg\{ w(G)-\dfrac{\varepsilon}{2}\,,0\,\bigg\}\Bigg\}\\
&\qquad\quad\qquad\qquad\qquad\qquad\qquad\qquad\qquad\qquad\mbox{and}\\
&\quad\qquad
w''(G)=\underset{G_*,G*\in\mathcal{T}\setminus\{\varnothing\}}{\min_{{G_*}\subseteq G\subseteq{G^*}}}\dfrac{\Big(w({G_*})+\varepsilon/2\Big)+\Big(w({G^*})+\varepsilon/2\Big)}{2}
\leq w(G)+\dfrac{\varepsilon}{2}\,.
}
Utilizing these two inequalities in \eq{700}, we get \eq{369}. This completes the proof.
\end{proof}
The present research also opens up several possibilities for further exploration. For example, the core concept of this paper can also be implemented on appropriate lattice structures. Besides one can replace the fixed error term $\varepsilon$ with a properly designed error function to obtain more generalized and interesting results.

\end{document}